\documentclass{ifacconf}

\usepackage{graphicx}      
\usepackage{natbib}        
\usepackage{epstopdf}
\DeclareGraphicsExtensions{.eps}
\usepackage{amsmath,amsfonts,amssymb}
\usepackage{dsfont}
\usepackage{color}

\newtheorem{definition}[thm]{Definition}

\newcommand{\Ical}{\mathcal{I}}
\newcommand{\Ncal}{\mathcal{N}}
\newcommand{\N}{\mathds{N}}
\newcommand{\R}{\mathds{R}}

\newcommand{\T}{\mathrm{T}}
\newcommand{\ones}{\mathds{1}}

\begin{document}
\begin{frontmatter}

\title{Scalable robustness of interconnected systems subject to structural changes}

\author[First]{Steffi Knorn} and
\author[Second]{Bart Besselink} 

\address[First]{Institute of Automation Engineering, Otto-von-Guericke University Magdeburg, Germany, Email: steffi.knorn@ovgu.de}
\address[Second]{Bernoulli Institute for Mathematics, Computer Science and Artificial Intelligence and the Jan C.\ Willems Centre for Systems and Control, University of Groningen, The Netherlands}

\begin{abstract}  
This paper studies the robustness of large-scale interconnected systems with respect to external disturbances, focussing on their scalability properties. Specifically, a notion of scalability is introduced that asks for these robustness properties to remain unchanged under a structural change of the system, such as the addition/removal of a subsystem or a change in the interconnection structure. Both necessary and sufficient conditions, in terms of the interconnection structure and edge weights, are given under which elementary structural changes are scalable. The results are illustrated through a simple example.
\end{abstract}

\end{frontmatter}

\section{Introduction}%
Technical systems have been evolving and increasing both in size and complexity over the last decades. This trend is accelerated by connecting a growing number of subsystems by means of communication, i.e., information exchange, or physical interconnections. Such large-scale interconnected systems can indeed be found in a variety of applications, e.g., formations of automatic vehicles, electricity networks including smart grids, sensor networks, and traffic systems.

The stability analysis and control of these interconnected systems has a long history, e.g., \cite{moylan_1978}, \smash{\cite{book_siljak_1978}}, whereas also consensus properties have been studied extensively, see \cite{fax_2004,olfati-saber_2007}, and \cite{knorn_2016b} for examples.

In this paper, however, we are interested in the robustness of interconnected systems with respect to external disturbances. We will study how these robustness properties change under structural changes of the  system, i.e., the addition or removal of subsystems or changes to the interconnection structure. This is motivated by, first, the observation that external disturbances on interconnected systems propagate through the network and can therefore cause undesirable effects ranging from loss of performance to cascaded failures. Second, modern large-scale engineering systems are generally subject to change, where the introduction of new generators or power transmission lines in electricity networks is an example. Such structural changes have an influence on the properties of the network.

The effect of external disturbances (sometimes called network performance) has been studied for \emph{static} interconnected systems in, e.g., \cite{zelazo_2011}, \cite{siami_2016}, and \cite{lovisari_2013}, whereas \cite{siami_2018} considers the effects of adding edges in the interconnection structure. An alternative perspective is given in \cite{bamieh_2012} (see also \cite{tegling_2019}), where network performance measures are analysed as a function of the network size, i.e., the number of subsystems, in a property known as \emph{coherence}. This is thus \emph{scalability} analysis of network performance.

An important example of a scalable performance notion is \emph{string stability} characterizing the amplification of disturbances in vehicle-following systems, e.g., platoons, see \cite{alam_2015b}. Even though various definitions of string stability exist, see \cite{fenton_1968,swaroop_1996,seiler_2004,barooah_2005,middleton_2010} and \cite{knorn_2014}, these notions have in common that they ask for an upper bound on a platoon performance measure \emph{independently} of the platoon size.
For extensions of string stability relying on the existence of a uniform performance bound for more general interconnection structures see \cite{knorn_2016} and \cite{besselink_2018}.

Even though string stability notions and coherence allow for characterising the effect of external disturbances for varying network size, they have in common that are only applicable for highly homogeneous systems. Also, more importantly, only highly structured interconnection topologies are considered.

This paper addresses these two limitations by, first, considering heterogeneous subsystem dynamics and, second, allowing for arbitrary changes in the interconnection structure, leading to the following contributions.

First, we introduce the notion of \emph{scalability} of structural changes, which are changes of the interconnected system due to removing or adding subsystems or changing the interconnection structure. Scalability requires that network performance is not decreased despite this structural change, where performance is understood as a bounded gain from external disturbances to state deviations.

Second, we analyse the scalability of elementary structural changes, namely the addition/removal of a subsystem or the addition/removal of an edge in the interconnection structure. It is shown that the removal of subsystems and edges is always scalable and both necessary and sufficient conditions are given for scalability of the addition of subsystems and edges. The latter results rely on relating robustness properties to walks in an associated graph and have a connection to small-gain conditions.

This paper is organised as follows. The problem statement is found in Sec.~\ref{sec_problem}. After characterising our notion of $\gamma$-robustness of large-scale systems in Sec.~\ref{sec:robustness}, results on scalable changes are presented in Sec.~\ref{sec_changes}. The paper closes with an example in Sec.~\ref{sec_example} and conclusions in Sec.~\ref{sec_conclusions}.

\emph{Notation}: The sets of natural numbers and real numbers are denoted $\N$ and $\R$. We write $\mathcal{I}_N = \{1,2,\ldots,N\}$. The all-ones vector is denoted $\mathds{1}$ and $e_i$ is the $i$th column of the identity matrix. For a matrix $A$, we denote its element in row $i$ and column $j$ by $(A)_{ij}$. A diagonal matrix with diagonal entries $a_i$ is denoted $\mathrm{diag}(a_1, a_2,\dots, a_N)$. A matrix $A$ is said to be Metzler if all its off-diagonal elements are nonnegative; $B$ is an M-matrix if $B = -A$, where $A$ is both Metzler and Hurwitz. We denote $\|d\|_\infty = \sup_t |d(t)|$. Vector inequalities are to be understood elementwise.

\section{Problem statement}\label{sec_problem}
Consider the linear large-scale, interconnected system comprised of $N$ scalar subsystems of the form
\begin{equation}\label{eqn_sysi}
	\Sigma_i: \dot{x}_i = -a_i x_i +\sum_{j\in \mathcal N_i} m_{ij}x_j +d_i,
	\end{equation}
with state $x_i\in\R$ and external disturbance $d_i\in\R$. Furthermore, $a_i>0$ for all $i$ are the self-feedback parameters, $m_{ij}>0$ for all $i,j$ capture the connecting edges, and $\mathcal N_i$ captures all (in-)neighbours of node $i$, which have edges pointing towards $i$ (not including $i$ itself). Then, the entire system can be described compactly by
\begin{equation}\label{eqn_sys}
	\Sigma: \dot{x} = -(A-M)x +d,
	\end{equation}
where $A=\mathrm{diag}(a_1, a_2,\dots, a_N)\in\R^{N\times N}$ and $M\in\R^{N\times N}$ is a matrix with nonnegative off-diagonal entries $(M)_{ij}=m_{ij}$ if $j\in\mathcal N_i$, $(M)_{ij}=0$ otherwise, and $(M)_{ii}=0$.

We are interested in structural changes to the system~\eqref{eqn_sys}, i.e., the addition/removal of subsystems, and/or connections. Subsystems and connections will also be referred to as nodes and edges, respectively, interpreting the interconnection structure as a graph. Such changes lead to the new system
\begin{equation}\label{eqn_sys_changed}
	\bar\Sigma: \dot{\bar{x}} = -(\bar A-\bar M)\bar{x} + \bar{d},
	\end{equation}
with potentially different matrices $\bar A$ or $\bar M$ depending on the particular change under consideration. Specifically, we consider four \emph{elementary} structural changes.
	
\textit{1. Removal of a node (without any attached edges).} Assume, without loss of generality, that the node to be removed is labeled by $N$. This leads to the change
\begin{align}
A = \begin{bmatrix} \bar{A} & 0 \\ 0 & a_{N} \end{bmatrix}, \;
M = \begin{bmatrix} \bar{M} & 0 \\ 0 & 0 \end{bmatrix} \quad\mapsto\quad
\bar{A},\; \bar{M}.\label{eqn_change_removenode}
\end{align}

\textit{2. Removal of an edge.} Let $(i,j)$ with $j\in\Ncal_i$ be an edge with weight $m_{ij}$, such that its removal can be written as
\begin{align}
A,\; M \quad\mapsto\quad \bar{A} = A, \; \bar{M} = M - m_{ij}e_ie_j^{\T}. \label{eqn_change_removeedge}
\end{align}

\textit{3. Addition of a node (without any attached edges).} The addition of a node, which without loss of generality is labelled $N+1$, leads to the change
\begin{align}
A,\; M \quad\mapsto\quad
\bar{A} = \begin{bmatrix} A &0 \\ 0 & a_{N+1} \end{bmatrix}, \;
\bar{M} = \begin{bmatrix} M &0 \\ 0 & 0 \end{bmatrix}.\label{eqn_change_addnode}
\end{align}

\textit{4. Addition of an edge.} Assume that $(M)_{ij} = 0$. Then, the addition of the edge $(i,j)$ with weight $m_{ij}>0$ yields
\begin{align}
A,\; M \quad\mapsto\quad \bar{A} = A, \; \bar{M} = M + m_{ij}e_ie_j^{\T}.\label{eqn_change_addedge}
\end{align}

It is important to note that any structural change to \eqref{eqn_sys} can be expressed as a sequence of elementary structural changes. E.g., changing the dynamics or parameters of existing nodes or edges can be understood as first deleting the node or edge and then adding a new node or edge, respectively, with the new, desired properties.

In this paper, we study the effect of structural changes on system properties, namely, the robustness of system \eqref{eqn_sys} with respect to the external disturbances $d$:
\begin{definition}\label{def_gammarobust}
A system \eqref{eqn_sys} is said to be $\gamma$-robust (with $\gamma >0$) if it is asymptotically stable and
\begin{align}
\max_{i\in\Ical_N} |x_i(t)| \leq \gamma \max_{i\in\Ical_N} \|d_i\|_{\infty}
\label{eqn_def_gammarobust}
\end{align}
for all $t\geq0$ and all trajectories satisfying $x_i(0) = 0$.
\end{definition}
Thus, $\gamma$-robustness amounts to a bounded $L_{\infty}$ gain from the external disturbances to the state deviations. As such, it gives a bound on the worst-case perturbation $x_i$ in terms of the largest disturbance $d_i$ and can be regarded as a notion of performance. In many applications, it is undesirable for a structural change to lead to decreased system performance, which motivates the following definition.

\begin{definition}\label{def_scalablechange}
A structural change to the system \eqref{eqn_sys} is said to be $\gamma$-scalable if both $\Sigma$ in \eqref{eqn_sys} and $\bar{\Sigma}$ in \eqref{eqn_sys_changed} are $\gamma$-robust. We simply call a structural change scalable if, for all $\gamma$ for which $\Sigma$ is $\gamma$-robust, we have that $\bar{\Sigma}$ is $\gamma$-robust.
\end{definition}
Def.~\ref{def_scalablechange} above asks for robustness properties of the system (in terms of $\gamma$-robustness) to be preserved after structural changes. Thus, scalability implies the weaker notion of $\gamma$-scalability, which merely requires $\gamma$-robustness for the systems before and after any structural change. This allows $\bar{\Sigma}$ to be less robust with respect to disturbances than $\Sigma$ as long as the same pre-specified bound $\gamma$ holds.
\begin{rem}
The notion of ($\gamma$-)scalability in Def.~\ref{def_scalablechange} allows for considering large-scale systems with heterogeneous subsystems and arbitrary interconnection structures. Consequently, it can be regarded as a generalization of scalable robustness notions such as string stability (e.g., \cite{ploeg_2014,besselink_2017}) and scalable input-to-state stability (sISS) (see \cite{besselink_2018}). In these cases, it is required that there exists an upper bound on the gain from the largest disturbance to the worst-case perturbation, \emph{independently} of the system size.
In the terminology of Def.~\ref{def_scalablechange}, these notions thus ask for the preservation of robustness properties for \emph{sequences} of (elementary) changes. In contrast to Def.~\ref{def_scalablechange}, string stability and sISS heavily rely on uniformity in the subsystem dynamics and interconnection structure.
\end{rem}
In the remainder of this paper, we will investigate if (or, under which conditions) the elementary structural changes listed above are scalable according to Def.~\ref{def_scalablechange}.

\section{$\gamma$-robustness of large-scale systems}\label{sec:robustness}
Before studying the scalability of structural changes, the notion of $\gamma$-robustness in Def.~\ref{def_gammarobust} is characterized in this section. Here, we will frequently exploit the following necessary and sufficient condition for $\gamma$-robustness.
\begin{lem}\label{lem_gammarobust}
A system \eqref{eqn_sys} is $\gamma$-robust if and only if there exists $v\in\R^N$ such that $v>0$ and
\begin{align}
-(A-M)v + \ones \leq 0, \qquad v \leq \gamma\ones
\label{eqn_lem_gammarobust}
\end{align}
\end{lem}
\begin{pf}
The result is closely related to \cite[Prop.~4]{rantzer_2015} and \cite[Lem.~2]{briat_2013}, that consider a strict inequality in~\eqref{eqn_def_gammarobust}. For completeness, we give a full proof.

\textit{if)} By \eqref{eqn_lem_gammarobust}, we have that $-(A-M)v <0$, which implies that $-(A-M)$ is Hurwitz as $-(A-M)$ is a Metzler matrix. Take an input $d(\cdot)$ such that $\max_{i\in\Ical_N}\|d_i\|_{\infty} \leq 1$. Then,
\begin{align}
-(A-M)v \leq -\ones \leq d(t) \leq \ones \leq (A - M)v
\label{eqn_lem_gammarobust_proof_inputbound}
\end{align}
for all $t$. Following \cite[Prop.~4]{rantzer_2015}, consider three trajectories of the system \eqref{eqn_sys}:
\begin{itemize}
	\item[1.] $x_{\mathrm{max}}(\cdot)$ corresponding to initial condition $x(0) = v$ and input $d(t) = (A-M)v$;
	\item[2.] $x_{\mathrm{min}}(\cdot)$ corresponding to initial condition $x(0) = -v$ and input $d(t) = -(A-M)v$; and
	\item[3.] $x(\cdot)$ corresponding to initial condition $x(0) = 0$ and arbitrary input $d(\cdot)$ satisfying~\eqref{eqn_lem_gammarobust_proof_inputbound}.
\end{itemize}
For $x_{\mathrm{max}}$ and $x_{\mathrm{min}}$, the initial conditions are the equilibrium corresponding to the respective inputs, such that $x_{\mathrm{max}}(t) = -x_{\mathrm{min}}(t) = v$ for all $t\geq0$.
Now, given that $-(A-M)$ is Metzler, the system~\eqref{eqn_sys} is monotone and we have that, for all $t\geq0$, $-v = x_{\mathrm{min}}(t) \leq x(t) \leq x_{\mathrm{max}}(t) = v$, such that the condition $v\leq\gamma\ones$ implies \eqref{eqn_def_gammarobust}.

\textit{only if)} Consider the input $d(t) = \ones$, $t\geq0$ and define $v = \lim_{t\rightarrow\infty}\int_0^t e^{-(A-M)s}\ones \,\mathrm{d}s$. Note that the limit exists as a result of asymptotic stability (recall Def.~\ref{def_gammarobust}). In addition, as $A-M$ is an M-matrix, $v>0$, see \cite{book_berman_1994}. By \eqref{eqn_def_gammarobust} we immediately have $v\leq \gamma\ones$. Moreover, $v$ is an equilibrium of the system for constant input $d(t) = \ones$ and thus satisfies $0 = -(A-M)v + \ones$, implying~\eqref{eqn_lem_gammarobust}.\hfill$\qed$
\end{pf}

\begin{rem}
It follows from the proof of Lem.~\ref{lem_gammarobust} that, owing to monotonicity of the system \eqref{eqn_sys}, the condition~\eqref{eqn_def_gammarobust} in fact holds for all trajectories with initial conditions $|x_i(0)| \leq v_i$, where $v_i$ is the $i$th element in $v$.
\end{rem}

Whereas Lem.~\ref{lem_gammarobust} characterises $\gamma$-robustness for a given $\gamma>0$, we are generally interested in the smallest $\gamma$ that makes \eqref{eqn_sys} $\gamma$-robust. To characterise such $\gamma$, it is observed that the conditions of Lem.~\ref{lem_gammarobust} are equivalent with the matrix $-(A-M)$ being Hurwitz and the vector
\begin{align}
u = (A - M)^{-1}\ones
\label{eqn_vbar}
\end{align}
satisfying $u\leq\gamma\ones$ ($u>0$ is guaranteed by \eqref{eqn_vbar}). For any solution $v$ to \eqref{eqn_vbar}, we have $u \leq v$, such that the $\max_i u_i$ gives the smallest $\gamma$ for which \eqref{eqn_sys} is $\gamma$-robust.

To give an interpretation for the vector $u$ in \eqref{eqn_vbar}, denote by $\mathcal{G}(MA^{-1})$ the directed weighted graph characterised through the weighted adjacency matrix $MA^{-1}$. Thus, $\mathcal{G}(MA^{-1})$ inherits the interconnection structure from $M$ in \eqref{eqn_sys} but scales the weights of the edges pointing from a node $i$ with $a_i^{-1}$. Specifically, $\mathcal{G}$ has the vertex set $\mathcal{I}_N$ and the set of edges $\mathcal{E}$ satisfies $(i,j)\in\mathcal{E}$ if and only if $m_{ij} > 0$.

We recall that, for a directed weighted graph $\mathcal{G}$, a (directed) \emph{walk} (of length $k$) from node $i$ to $j$ is a sequence of nodes $(i_0,i_1,\ldots,i_k)$ such that $i_0 = i$, $i_k = j$ and $(i_{l+1},i_l)\in\mathcal{E}$. To such walk we associate a weight given by the product of the weights of all edges that are traversed, which for a walk in $\mathcal{G}(MA^{-1})$ amounts to $\prod_{l=0}^{k-1} \frac{m_{i_{l+1},i_l}}{a_{i_{l}}}$. A \emph{path} is a walk for which all nodes in the sequence are distinct. A path with $i = j$ is called a \emph{cycle}. We then obtain the following interpretation for $u$.

\begin{lem}\label{lem_weightedwalks}
Consider the system \eqref{eqn_sys}. Then, $-(A - M)$ is Hurwitz if and only if $\rho(MA^{-1}) < 1$. In this case, the element $u_j$ of \eqref{eqn_vbar} is such that $u_ja_j - 1$ equals the sum of all weighted walks in the graph $\mathcal{G}(MA^{-1})$ that end in~$j$.
\end{lem}
\begin{pf}
The first statement can be found in \cite[Thm.~30]{duan_2019_arxiv}. To prove the second, note that
\begin{align}
u = (A-M)^{-1}\ones = A^{-1}(I - MA^{-1})^{-1}\ones.
\label{eqn_lem_weightedwalks_proof_step1}
\end{align}
As $\rho(MA^{-1}) < 1$, we have (e.g., \cite[Lem.~2.1]{book_berman_1994}) that
\begin{align}
(I - MA^{-1})^{-1} = I + \sum_{k=1}^{\infty} (MA^{-1})^k.
\label{eqn_lem_weightedwalks_proof_step2}
\end{align}
It is well-known that $((MA^{-1})^k)_{ji}$ for $k>0$ is the sum of the weighted walks of length $k$ from node $i$ to node $j$, which proves the result through \eqref{eqn_lem_weightedwalks_proof_step1}.\hfill$\qed$
\end{pf}
\begin{rem}
The result in Lem.~\ref{lem_weightedwalks} shows that two factors contribute to the size of entry $j$ in $u$.
\begin{itemize}
	\item The number of incoming edges, paths, or walks. Specifically, $u_j$ increases when node $j$ is influenced heavily by other nodes (characterised through the number of walks ending in $j$). Roughly speaking, such influence causes the disturbances acting on nodes other than $j$ to propagate through the network to~$j$.
	\item The weight of the self-feedback parameter $a_j$. Namely, $u_j$ decreases for increasing $a_j$, in which case subsystem $j$ has itself an increased robustness with respect to incoming disturbances (both directly through $d_j$ as those that have propagated through the network).
\end{itemize}
Thus, the entries in $u$ and hence the bound $\gamma$ in \eqref{eqn_def_gammarobust} can be kept small by limiting the number of directed paths and edges, specially ending in one particular node, and using sufficiently high self-feedback parameters.
\end{rem}

The interpretation of $u$ in terms of weighted walks in the graph $\mathcal{G}(MA^{-1})$ also allows for obtaining the following necessary condition for $\gamma$-robustness, expressed in terms of a small-gain type condition for cycles in $\mathcal{G}(MA^{-1})$.
\begin{lem}\label{lem_smallgain}
Consider the system \eqref{eqn_sys} and assume that $-(A-M)$ is Hurwitz. Then, \eqref{eqn_sys} is $\gamma$-robust only if for each cycle in $\mathcal{G}(MA^{-1})$, its weight $w$ satisfies
\begin{align}
\frac{1}{1 - w} \leq a_i\gamma
\label{eqm_lem_smallgain}
\end{align}
for all $i$ such that node $i$ is part of the cycle.
\end{lem}

\begin{pf}
Consider a cycle of length $k$ in $\mathcal{G}(MA^{-1})$, let $i$ be any node in this cycle, and denote by $w$ the weight associated to this cycle. Then, it necessarily holds that
\begin{align}
\left( (MA^{-1})^k \right)_{ii} \geq w,
\label{eqn_lem_smallgain_proof_step1}
\end{align}
as the cycle is amongst the walks of length $k$ from $i$ to $i$. Then, the use of \eqref{eqn_lem_weightedwalks_proof_step1} and \eqref{eqn_lem_weightedwalks_proof_step2} leads to $a_iu_i = e_i^{\T}(I - MA^{-1})\ones = 1 + e_i^{\T}\left( \sum_{l=1}^{\infty} (MA^{-1})^l \right)\ones$, whereas the observation that $MA^{-1}\geq0$ allows for showing the sequence of lower-bounds
\begin{align}
\nonumber a_iu_i &\geq 1 + \sum_{l=1}^{\infty} \left(\left(MA^{-1}\right)^l\right)_{ii} 
\geq 1 + \sum_{l=1}^{\infty} \left(\left(MA^{-1}\right)^{kl}\right)_{ii}, \\
&\geq 1 + \sum_{l=1}^{\infty} \left(\left( \left(MA^{-1}\right)^k \right)_{ii}\right)^l. \label{eqn_lem_smallgain_proof_step2}
\end{align}
Now, after recalling \eqref{eqn_lem_smallgain_proof_step1}, \eqref{eqn_lem_smallgain_proof_step2} yields $a_iu_i \geq 1 + \sum_{l=1}^{\infty} w^l = \frac{1}{1 - w}$, where the equality follows as the inequality implies that the sum converges. Then, noting that $u_i\leq\gamma$, we obtain the necessary condition \eqref{eqm_lem_smallgain}.\hfill$\qed$
\end{pf}
\begin{rem}
Condition \eqref{eqm_lem_smallgain} can be regarded as a small-gain condition. Stability and input-to-state stability of large-scale interconnected nonlinear systems have been studied extensively using small-gain conditions, see \cite{dashkovskiy_2007,dashkovskiy_2010c}. These results however differ from the result in Lem.~\ref{lem_smallgain} as the latter gives a small-gain result that is necessary for \emph{a given} robustness bound $\gamma$, whereas the former generally target stability or input-to-state stability (but without explicitly characterising the corresponding gain functions).
\end{rem}

\section{Scalable structural changes}\label{sec_changes}
In this section, scalability of \emph{elementary} structural changes according to Def.~\ref{def_scalablechange} is considered.

First, consider the removal of a node without any attached edges.
\begin{prop}\label{lem_delete_note}
Let the system $\Sigma$ in \eqref{eqn_sys} be $\gamma$-robust. Then, deleting a node without incoming or outgoing edges, characterised through \eqref{eqn_change_removenode}, is a scalable change.
	\end{prop}

\begin{pf}
Consider any $\gamma$ such that $\Sigma$ is $\gamma$-robust. By Lem.~\ref{lem_gammarobust}, there exists a vector $v\in\R^N$ satisfying $0<v\leq\gamma\ones$ such that \eqref{eqn_lem_gammarobust} holds. After partitioning $v$ as $v = [\begin{array}{cc} \bar{v}^\T & v_N \end{array}]^\T$, the condition \eqref{eqn_lem_gammarobust} can be written as
\begin{equation}
	-\begin{bmatrix} \bar A-\bar M & 0 \\ 0 & a_N\end{bmatrix} 
	\begin{bmatrix} \bar v \\ v_N\end{bmatrix} 
	= -\begin{bmatrix} (\bar A-\bar M)\bar v \\ a_N v_N\end{bmatrix} 
	\leq -\mathds 1.
	\end{equation}
As the inequality is element-wise, this implies $-(\bar A-\bar M)\bar v + \mathds 1 \leq 0$ with $0<\bar{v}\leq\gamma\ones$, i.e., $\bar{\Sigma}$ is $\gamma$-robust.\hfill$\qed$
\end{pf}

Similar to the removal of nodes, we will further show that removing any edge is a scalable change.
\begin{prop}\label{lem_delete_edge}
Let the system $\Sigma$ in \eqref{eqn_sys} be $\gamma$-robust. Then, removing any edge, as in~\eqref{eqn_change_removeedge}, is a scalable change.
\end{prop}
\begin{pf}
Let $v$ satisfy the conditions of Lem.~\ref{lem_gammarobust} for the $\gamma$-robust system $\Sigma$. Then, using \eqref{eqn_change_removeedge}, we obtain
\begin{align}
	\nonumber -(A-\bar M)v =& -(A-M)v -m_{ij}e_ie_j^\mathrm Tv\\
	 = -(A-M)v -m_{ij}v_je_i
	\leq& -(A-M)v \leq -\mathds 1,
	\end{align}
where $v_j>0$ is the $j$th entry of $v$ and the final inequality follows from $\gamma$-robustness of $\Sigma$. As $\gamma$ is arbitrary, the change is scalable.\hfill$\qed$
\end{pf}

Hence, removing nodes or edges is always guaranteed to be a scalable change. The opposite is not true in general, as will be shown below. However, the addition of a node, denoted $N+1$, without incoming or outgoing edges can be made a scalable change by appropriate choice of the self-feedback parameter $a_{N+1}$.

\begin{prop}\label{lem_adding_node}
Let the system $\Sigma$ in \eqref{eqn_sys} be $\gamma$-robust and consider the associated vector $u$ in \eqref{eqn_vbar}. Then, adding a node without incoming or outgoing edges, as in~\eqref{eqn_change_addnode}, is a scalable change if and only if its self-feedback parameter $a_{N+1}$ satisfies $a_{N+1} \geq (\max_{i\in\Ical_N} u_i)^{-1}$.
\end{prop}
\begin{pf}
Consider the equation
\begin{equation}
	-(\bar A-\bar M) \bar{u}
	=-\begin{bmatrix}  A-M & 0 \\ 0 & a_{N+1}\end{bmatrix} 
	\begin{bmatrix} u \\ u_{N+1}\end{bmatrix} 
	= -\mathds 1, \label{eqn_lem_adding_node_proof_step1}
	\end{equation}
where we have introduced the partitioning $\bar{u} = [\begin{array}{cc} u^\T & u_{N+1} \end{array}]^\T$ and exploited the observation that $u$ satisfies \eqref{eqn_vbar}. The latter is due to $\Sigma$ being $\gamma$-robust, for any $\gamma$ such that $u \leq \gamma\ones$. We would like to use \eqref{eqn_lem_adding_node_proof_step1} as the counterpart of~\eqref{eqn_vbar} for the updated system $\bar{\Sigma}$, which requires, first, that $u_{N+1} > 0$. This implies $a_{N+1} > 0$. As a result, $-(\bar{A}-\bar{M})$ is invertible and the vector $\bar{u}$ indeed solves \eqref{eqn_vbar} for $\bar{\Sigma}$. Then, for the change to be scalable, it is required that $u_{N+1} \leq \gamma$ for any $\gamma$ satisfying $u\leq\gamma\ones$, which implies that $u_{N+1} \leq \max_{i\in\Ical_N} u_i$. This proves necessity.
Sufficiency follows immediately from \eqref{eqn_lem_adding_node_proof_step1}.\hfill$\qed$
\end{pf}

The following result gives a necessary and sufficient condition for the addition of an edge to be scalable.
\begin{prop}\label{lem_adding_edge}
Let the system $\Sigma$ in \eqref{eqn_sys} be $\gamma$-robust and consider the associated vector $u$ in \eqref{eqn_vbar}. Then, adding an edge $(i,j)$ (i.e., from node $j$ to $i$) with weight $m_{ij} > 0$, characterised through~\eqref{eqn_change_addedge}, is a scalable change if and only if $-(A - M - m_{ij}e_ie_j^{\T})$ is Hurwitz and
\begin{align}
\frac{m_{ij}}{1 - m_{ij}e_j^{\T}(A - M)^{-1}e_i} (A-M)^{-1}e_iu_j \leq u_{\textrm{max}}\ones - u,
\label{eqn_lem_adding_edge}
\end{align}
where $u_{\textrm{max}} = \max_{k\in\Ical_N} u_k$.
\end{prop}
\begin{pf}
As the vector $u$ characterises the smallest $\gamma$ for which $\Sigma$ in \eqref{eqn_sys} is $\gamma$-robust, it is clear that the change is scalable if and only if $-(A - M - m_{ij}e_ie_j^{\T})$ is Hurwitz and $\bar{u} = (A - M - m_{ij}e_ie_j^{\T})^{-1}\ones \leq u_{\textrm{max}}\ones$. However, the Sherman-Morrison formula for matrix inverses (e.g., \cite{book_horn_2013}) gives
\begin{align}
\bar{u}
&= (A-M)^{-1}\ones + \frac{m_{ij}(A-M)^{-1}e_ie_j^{\T}(A-M)^{-1}}{1 - m_{ij}e_j^{\T}(A-M)^{-1}e_i} \ones, \nonumber\\
&= u + \frac{m_{ij}}{1 - m_{ij}e_j^{\T}(A-M)^{-1}e_i}(A-M)^{-1}e_iu_j,
\label{eqn_lem_adding_edge_proof_step2}
\end{align}
where \eqref{eqn_lem_adding_edge_proof_step2} is a result of the definition of $u$ in \eqref{eqn_vbar}. The result \eqref{eqn_lem_adding_edge} follows after requiring that $\bar{u} \leq u_{\textrm{max}}\ones$.\hfill$\qed$
\end{pf}

Although \eqref{eqn_lem_adding_edge} is not easily verified in practice for large-scale interconnected systems, it allows for an insightful interpretation. Namely, by noting that $-(A-M)$ is Metzler and Hurwitz, it follows that $(A-M)^{-1}$ has nonnegative entries (this can also be observed from Lem.~\ref{lem_weightedwalks}). Consequently, $u\leq\bar{u}$ and the addition of an edge can never lead to stronger robustness properties. In fact, the left-hand side of \eqref{eqn_lem_adding_edge} characterises the performance loss.

This performance loss is determined by two factors, which can be understood as follows. First, $e_j^\T(A-M)^{-1}e_i$ can be interpreted as the weighted sum of all walks from node $i$ to $j$. Recalling that we have added an edge from node $j$ to $i$, this additional edge introduces cycles with respect to those walks, after which the condition $m_{ij}e_j(A-M)^{-1}e_i < 1$ can be regarded as a small-gain like condition.

Second, the $k$th element in $(A-M)^{-1}e_i$ gives the sum of all weighted walks from node $i$ to node $k$ and thus characterises how the effect of disturbances that can be propagated through the added edge $(i,j)$ are distributed through the network. This depends on how much node $j$ was affected by disturbances in the first place, given by $u_j$.

The above observations lead to a few corollaries describing cases in which the addition of an edge cannot be scalable.
\begin{cor}
Let \eqref{eqn_sys} be $\gamma$-robust, consider the associated vector $u$ in \eqref{eqn_vbar}, and let $k$ be any index such that $u_k = \max_{i\in\Ical_N} u_i$. Consider adding the edge $(i,j)$. If there exists a directed path from $i$ to $k$, the change is not scalable.
\end{cor}

Next, we provide some sufficient conditions under which the addition of an edge is scalable.
\begin{prop}\label{lem_adding_edge_sufficient}
Let~\eqref{eqn_sys} be $\gamma$-robust and consider an associated vector $v$ satisfying \eqref{eqn_lem_gammarobust}. Then, the addition of the edge $(i,j)$ with weight $m_{ij}>0$ is scalable if
\begin{equation}
	a_i \geq \frac{\sum_{k\in{\mathcal N}_i} m_{ik}v_k + m_{ij}v_j +1}{v_i}.
	\label{eqn_lem_adding_edge_sufficient}
	\end{equation}
\end{prop}
\begin{pf}
Consider the change \eqref{eqn_change_addedge} and note that, by Lem.~\ref{lem_gammarobust}, this change is scalable if
\begin{align}
-(A - \bar{M})v + \ones \leq 0, \label{eqn_lem_adding_edge_sufficient_proof_step1}
\end{align}
where we have used the vector $v>0$ associated to $\Sigma$. As $\bar{M} = M + m_{ij}e_ie_j^{\T}$, we have that $e_l^{\T}(A - \bar{M}) = e_l^{\T}(A - M)$, for all $l\neq i$, such that the element-wise inequality \eqref{eqn_lem_adding_edge_sufficient_proof_step1} is guaranteed to hold for all rows $l\neq i$. For row $l=i$, \eqref{eqn_lem_adding_edge_sufficient_proof_step1} reads
\begin{align}
-a_iv_i + \sum_{k\in\Ncal_i} m_{ik}v_k + m_{ij}v_j + 1 \leq 0,
\label{eqn_lem_adding_edge_sufficient_proof_step2}
\end{align}
which is guaranteed to hold under the assumption \eqref{eqn_lem_adding_edge_sufficient}. We remark that the term $m_{ij}v_j$ in \eqref{eqn_lem_adding_edge_sufficient_proof_step2} captures the effect of the added edge.\hfill$\qed$
\end{pf}

The result in Prop.~\ref{lem_adding_edge_sufficient} is indeed only sufficient as it is based on choosing the same $v$ for the system before and after the structural change. {Nonetheless, adding an edge might be scalable despite updating $v$ (to $\bar{v}$) as long as $\max_i v_i = \max_i \bar{v}_i$. The utility of Prop.~\ref{lem_adding_edge_sufficient} is clear after noting that \eqref{eqn_lem_adding_edge_sufficient} can be verified locally, assuming that each subsystem $i$ has access to its relevant component $v_i$ of the vector $v$ characterising $\gamma$-robustness of the original system.}
\begin{rem}\label{rem_adjustselffeedback}
Note that, in practice, \eqref{eqn_lem_adding_edge_sufficient} might not be satisfied for all changes. However, one might choose to adapt the self-feedback parameter of node $i$ to a larger value $\bar{a}_i \geq a_i$ to render the change scalable. {Specifically, the change is scalable if $\bar{a}_i$ is chosen such that
\begin{align}
- \bar{a}_i v_i + \sum_{l\in\mathcal N_i} m_{il}v_l + m_{ij}v_j + 1\leq 0,
\end{align}
as follows immediately from \eqref{eqn_lem_adding_edge_sufficient}.} For instance, adding any edge $m_{ij}$ is scalable if setting $\bar a_i = a_i + m_{ij}\frac{v_j}{v_i}$.
\end{rem}

\begin{cor}
Let $\Sigma$ in \eqref{eqn_sys} be $\gamma$-robust, consider the associated vector $u$ in \eqref{eqn_vbar}. Then, the addition of the edge $(i,k)$ results in vector $\bar u$ for which $\bar u_i \geq u_i +\frac{m_{ik}}{a_i}u_k$.
\end{cor}
\begin{pf}
The result follows from the discussions above.
\end{pf}
\begin{rem}
The sufficient condition in Prop.~\ref{lem_adding_edge_sufficient} sheds an interesting light on the effects of adding an edge and supports the intuition already identified above. Namely, the effects of adding an edge (scaled with the weight of the edge $m_{ij}$) are larger in case the edge goes from a node heavily affected by disturbances to a node less affected by disturbances, i.e., when $v_j>v_i$.
\end{rem}

Even though it can in principle be verified in a decentralised manner, the condition for a scalable change in Prop.~\ref{lem_adding_edge_sufficient} could be difficult to check in practice as all $v_i$ have to be known. However, in some cases, the condition for adding an edge being scalable can be simplified significantly as will be shown below.
\begin{cor}\label{cor_adding_edge_sufficient}
Consider the system $\Sigma$ in \eqref{eqn_sys} and assume that the matrix $A-M$ is strictly diagonally dominant, i.e., $a_i>\sum_{j\in\Ncal_i}m_{ij}$ for all $i$. Then, there exists a vector $v=\gamma\mathds{1}$ satisfying condition \eqref{eqn_lem_gammarobust}. Further, the addition of the edge $(i,j)$ with weight $m_{ij}>0$ is scalable if
	\begin{equation}
		a_i \geq 
		\sum_{k\in{\mathcal N}_i} m_{ik} + m_{ij} +\frac{1}{\gamma}.
		\label{eqn_cor_adding_edge_sufficient}
		\end{equation}
\end{cor}
\begin{pf}
The existence of a vector $v = \gamma\ones$ satisfying \eqref{eqn_lem_gammarobust} follows immediately after noting that strict diagonal dominance implies the existence of a scalar $\gamma>0$ such that $a_i > \frac{1}{\gamma} + \sum_{j\in\Ncal_i} m_{ij}$, for all $i\in\mathcal{V}$, which can be rewritten to \eqref{eqn_lem_gammarobust}. The condition \eqref{eqn_cor_adding_edge_sufficient} is merely a restatement of \eqref{eqn_lem_adding_edge_sufficient}  for $v = \gamma\ones$.\hfill$\qed$
\end{pf}
Condition \eqref{eqn_cor_adding_edge_sufficient} can be verified in a decentralised manner assuming that the desired performance level $\gamma$ is known. We also note that the observations in Rem.~\ref{rem_adjustselffeedback} still hold.

The results above hence show an interesting insight into the structure of networks and the effects of adding edges. Any additional edge in a system, which is not ``compensated'' for by appropriate adjustment of the corresponding local self-feedback parameter, will have a negative impact on the network in the sense that the corresponding norm cannot decrease. (Unless of course in cases where by choosing bounds and corresponding vectors $v$ much larger than needed a sufficiently large margin can be guaranteed.) Further, $v_i$ and its increase capture the accumulated, weighted influence of the potential noise and disturbances acting on any node in the network for which a directed path exists to $i$ such that noise or disturbances acting on such a node will also affect $i$.

\section{Example}\label{sec_example}
Consider a network with three nodes and three edges such that the system is described by
\begin{equation}
	\dot{x} = \begin{bmatrix} -1 & 0 & 0 \\
							   1 & -1 & 0 \\
							   1 & 1 & -1\end{bmatrix}
	x + d,
	\end{equation}
which is $\gamma$-robust with the associated vector $u=\left(1,2,4\right)^\mathrm T$ such that $\gamma = 4$.

Adding a node is hence scalable as long as the self-feedback parameter is chosen larger or equal to $1/\gamma = 1/4$. Indeed, this would lead to system
\begin{equation}
	\dot{x} = \begin{bmatrix} -1 & 0 & 0 & 0 \\
							   1 & -1 & 0 & 0 \\
							   1 & 1 & -1 & 0\\
							   0 & 0 & 0 & -\tfrac{1}{4} \end{bmatrix}
	x + d, \quad\text{with }
	u = \begin{bmatrix} 1 \\ 2 \\ 4 \\ 4\end{bmatrix},
	\end{equation}
yielding the same upper bound $\gamma = 4$.

If an edge should also be added, consider for instance an edge from node 2 to 4. Since due to our choice $a_4 = \frac{1}{4}$ we have $v_4=\gamma$, any weight $m_{4,2}>0$ of any edge pointing towards node 4 or creating an additional path ending in node 4 will lead to an unscalable change. For instance, setting $m_{4,2}=0.1$ increases $v_4$ and hence also $\gamma$ to 4.8 which is larger than the previous bound 4.

Using Rem.~\ref{rem_adjustselffeedback}, the change of adding the edge from node 2 to 4 with weight $m_{4,2}=0.1$ can in fact be made scalable by adjusting the self-feedback parameter of node 4 to $\bar a_4 = \frac{1}{4}+0.1\frac{2}{4} = 0.3$.

\section{Conclusions}\label{sec_conclusions}
The robustness of large-scale interconnected systems with respect to disturbances is studied in this paper, focussing on the scalability of robustness properties with respect to structural changes to the system.

Future work will focus on the decentralised verification of scalability with respect to structural changes and the extension of the results in this paper to nonlinear systems.

\bibliography{references_arxiv_clean}

\end{document}